\documentclass[12pt]{amsart}
\usepackage{amssymb,mathrsfs,longtable}
\usepackage[matrix,arrow,curve]{xy}
\usepackage{setspace}
\usepackage{multirow}
\usepackage{textcomp}
\usepackage{hyperref}

\usepackage[dvipsnames]{xcolor}
\usepackage{graphicx}

\usepackage[normalem]{ulem}
\usepackage{cancel}

\newcounter{cequation}[section]

\newtheorem{theorem}[cequation]{Theorem}
\newtheorem*{theorem*}{Theorem}
\newtheorem{lemma}[cequation]{Lemma}

\theoremstyle{definition}
\newtheorem{example}[cequation]{Example}
\newtheorem{definition}[cequation]{Definition}
\newtheorem*{definition*}{Definition}

\newtheorem*{notation*}{Notation}

\theoremstyle{remark}
\newtheorem{remark}[cequation]{Remark}

\makeatletter\@addtoreset{equation}{section}
\makeatletter\@addtoreset{section}{part}

\makeatother

\def \O {\mathcal{O}}
\def \CC {\mathbb{C}}

\def \PP {\mathbb{P}}

\def \ge {\geqslant}
\def \le {\leqslant}

\def \Sing {\mathrm{Sing}\,}
\newcommand{\codim}{\mathrm{codim}\,}
\def \wt {\mathrm{wt}}
\newcommand{\Proj}{\mathrm{Proj}\,}
\newcommand{\Aff}{{\mathbb A}}
\newcommand{\rk}{\mathrm{rk}\,}

\title{Well formedness vs weak well formedness}

\author{Victor Przyjalkowski} % and Constantin Shramov}

\address{\emph{Victor Przyjalkowski}
\newline
\textnormal{Steklov Mathematical Institute of RAS, 8 Gubkina street, Moscow 119991, Russia.
}
\newline
\textnormal{\texttt{victorprz@mi-ras.ru, victorprz@gmail.com}}}
\thanks{This work is supported by the Russian Science Foundation under grant \textnumero 19-11-00164.}

\begin{document}

\begin{abstract}
In the literature there are two definitions of well formed varieties in weighted projective spaces.
According to the first one, well formed variety is the one whose intersection with the singular
locus of the ambient weighted projective space has codimension at least two,
while, according to the second one,  well formed variety is the one who does not contain in codimension one
a singular stratum of the ambient weighted projective space. We show that these two definitions
indeed differ, and show that they coincide for quasi-smooth weighted complete intersections of dimension
at least $3$.
\end{abstract}

\maketitle
%\tableofcontents

\section{Motivation}

One of the most straightforward examples of projective algebraic varieties are complete intersections
in projective spaces.
%, which, in some sense, can be viewed as ones given by general enough homogenous polynomials
%of fixed degrees.
%In particular, general complete intersections are smooth.
The generalization of usual complete intersections are weighted ones,
that are zeroes of weighted homogenous polynomials that form a regular sequence. However, due to
numerical properties of weights and degrees determining the complete intersection,
there can be not so much parameters to vary, so they can have some pathologies.
In particular, they can be singular for any choice of parameters. There are two reasons for a variety in weighted projective space
to be singular: singularities of the ambient weighted projective space and ``its own'' singularities.
If the variety does not have ``its own'' singularities, it is called \emph{quasi-smooth} and under some assumption behaves, from a lot of points
of view, as a smooth one (or, alternatively, it can be viewed as a smooth orbifold).
``Expected behavior'' with respect to singularities of the ambient weighted projective space is called \emph{well formedness}.
Informally it means that a subvariety in a weighted projective space is actually singular in the singular
points of the weighted projective space --- as it is natural to expect.

The standard classical references on weighted complete intersections are papers
of Dolgachev~\cite{Do82}, Iano-Fletcher~\cite{IF00}, and Dimca~\cite{Di86}. Well formedness condition is missed in~\cite{Do82}.
It was introduced in~\cite{Di86} and in~\cite{IF00}, and is used widely
for subvarieties of weighted projective spaces. The problem is that the two
definitions are close but different, which makes a mess in using them.
To fix this mess we show that these two definitions indeed differ,
but they coincide in an important ``good'' case. %of {quasi-smooth} weighted complete intersections of dimension at least $3$.
We call the condition from~\cite{Di86} well formedness (see Definition~\ref{definition:well formed wci}),
while the one from~\cite{IF00} we refer to as weak well formedness (see Definition~\ref{definition: weak well formed}).
Well formed varieties are weakly well formed. We give examples
(Examples~\ref{example: smooth not well formed} and~\ref{example: weakly not well formed})
when the converse is not true. However we show (Theorem~\ref{theorem: weakly wf is wf})
that it is true for the quasi-smooth weighted complete intersections of dimension at least $3$.
This justifies using results from~\cite{IF00} and the other papers
in this case.

\medskip

{\bf Acknowledgements.} This note is a byproduct of the work on the book ``Weighted complete intersections''
by the author and C.\,Shramov. I am grateful to him for the discussions without which
the note would not be written.

\section{Results}
Let $S=\CC[x_0,\ldots,x_N]$ be a graded polynomial ring given by assigning the weights
$\wt (x_i)=a_i$, and let $\PP=\PP(a_0,\ldots,a_N)=\Proj S$ be a weighted projective space. % of dimension $N$.
The weights of a weighted projective space are not invariant.
For instance, scaling grading gives the isomorphism
$$
\PP(aa_0,\ldots,aa_N)\simeq \PP(a_0,\ldots,a_N).
$$
The collection of weights may be further simplified.

\begin{definition}
\label{definition:P well formed}
We say that $\PP$ is \emph{well formed} if $\mathrm{gcd}(a_0,\ldots, a_{i-1},a_{i+1},\ldots, a_{N})=1$
for all $i=1,\ldots, N$.
\end{definition}
Note that the well formedness assumption is not restrictive.
Indeed, for any positive integer $a$ consider the Veronese embedding given by the inclusion
$$
\bigoplus_{i=0}^\infty S_{ai}\hookrightarrow \bigoplus_{i=0}^\infty S_{i},
$$
where $S_j$ is the $j$-th graded component of $S$. It gives the isomorphism
$$
\PP(a_0,a\cdot a_1,\ldots,a\cdot a_N)\cong\PP(a_0,a_1,\ldots,a_N),
$$
where $a_0$ and $a$ are coprime %(this isomorphism is again given by the Veronese map).,
%because by~
(see details in~\cite[1.3.1]{Do82} and~\cite[Lemma 5.7]{IF00}).
Thus every weighted projective space is isomorphic to a well formed one.
This shows that the well formedness condition for $\PP$ is not geometric, it is a convention
of numerical expression of $\PP$ (we consider weighted projective spaces as varieties,
not as orbifolds).

In the opposite to the usual projective spaces, weighted ones are singular unless all weights are equal to $1$.

\begin{lemma}[{see, for instance,~\cite[5.15]{IF00}}]
\label{lemma:singularities-of-P}
Let $\PP$ be a well formed weighted projective space.
For every subset
$$
J\subset \{0,\ldots,N\}
$$
such that the greatest common divisor $a_J$ of the weights~$a_j$
for~\mbox{$j\in J$} is greater than~$1$, denote
$$
\Lambda_J=\left\{(x_0:\ldots:x_N) \mid x_j=0 \text{\ for all\ } j\notin J\right\}.
$$
Then the singular locus of $\PP$ is a union of the strata $\Lambda_J$.
\end{lemma}

The weighted projective space comes with the natural projection
$$
p\colon \mathbb A^{N+1}\setminus \{0\}\to \PP
$$
given by the $\CC^*$-action.
The preimage $C^*_X=p^{-1}(X)$ of a subvariety~\mbox{$X\subset \PP$} is called the \emph{punctured affine cone of $X$}.
Its closure $C_X$ in $\mathbb A^{N+1}$ is called the \emph{affine cone of $X$}.
Due to the $\CC^*$-action, $C^*_X$ cannot have isolated singular points.

An important class of subvarieties in weighted projective spaces is formed by quasi-smooth ones.

\begin{definition}\label{definition: quasi-smoothness}
A subvariety $X$ of a weighted projective space $\PP$
is called \emph{quasi-smooth} if the punctured affine cone~$C^*_X$ is smooth.
\end{definition}

\begin{definition}
A subvariety $X\subset\PP$ of codimension $k\le N$ is called a \emph{weighted complete
intersection of multidegree $(d_1,\ldots,d_k)$} if it is given by an ideal $I=(f_1,\ldots,f_k)$ generated
by a regular sequence of homogenous polynomials $f_1,\ldots,f_k\in S$ of degrees $d_1,\ldots,d_k$, respectively.
In the other words, $X=\Proj S/I\subset \PP$.
\end{definition}

The following example shows that weighted complete intersection in a singular weighted
projective space can be smooth even if it passes through the singular locus of the ambient space.

\begin{example}\label{example:line-on-a-cone-not-WF}
Let $\PP=\PP(1,1,2)$, so that $\PP$ is a quadratic cone in $\PP^3$.
Let $X$ be a hypersurface of weighted degree $1$ in $\PP$, i.e. a line on this cone.
Then $X$ passes through the singular point of $\PP$, that is the origin of the
cone, but, obviously, $X$ is smooth.
\end{example}

To justify the expectation that complete intersection should be singular in
the singular locus of an ambient weighted projective space,
it was introduced the notion of \emph{well formedness}.

%The well formedness condition on $X$ has more challenging history.
%It was missed, but often implicitly used in~\cite{Do82}.
%Then it was introduced in~\cite{Di86} under the name of being in general position.
%Later it was introduced in~\cite{IF00}. However these definitions differ,
%which is confusing for using them. To stress this and fix the mess in using the two
%different notions in proofs we show that these two definitions indeed differ,
%but they coincide in the important case of \emph{quasi-smooth} varieties.

\begin{definition}[{\cite[Definition 1]{Di86}}]
\label{definition:well formed wci}
A subvariety $X$ of a weighted projective space $\PP$ is said to be \emph{well formed}
if~$\PP$ is well formed and
$$
\dim X-\dim \left( X\cap\Sing \PP \right)\ge 2,
$$
where the dimension of the empty set is defined to be $-1$.
\end{definition}

The line on a quadratic cone from Example~\ref{example:line-on-a-cone-not-WF} is not well formed variety. %The example of not well-formed complete intersection is given in Example~\ref{example:line-on-a-cone-not-WF};
Non-well formedness there is the reason of the unexpected smoothness. The similar phenomena appear in any dimension.

\begin{example}\label{example:line-on-quadratic-cone-generalized}
Generalizing Example~\ref{example:line-on-a-cone-not-WF}, let~\mbox{$\PP=\PP(1,1,2^{N-1})$}
be a weighted projective space with weighted homogeneous coordinates~\mbox{$x_0,\ldots,x_N$}, where $N\ge 2$.
Let~$X$ be the hypersurface of weighted degree~$1$ defined by equation $x_0=0$.
Then~$X$ passes through the singular locus of~$\PP$ defined by equations $x_0=x_1=0$ by Lemma~\ref{lemma:singularities-of-P};
moreover, $\Sing\PP$ has codimension $1$ in $X$, so that $X$ is not well formed. Note that $X\simeq \PP(1,2^{N-1})\simeq \PP^{N-1}$
and, in particular, is smooth.
\end{example}

In~\cite[6.9]{IF00} another notion of well formedness is introduced
(and consistently used, see e.g. the proof of \cite[Theorem~6.17]{IF00}).
We call it \emph{weak well formedness}.

\begin{definition}
\label{definition: weak well formed}
The variety $X\subset \PP$ is called \emph{weakly well formed} if $\PP$ is well formed and~$X$ contains
no singular linear strata
of $\PP$ whose codimension in $X$ equals~$1$.
\end{definition}

Obviously, every well formed variety is weakly well formed.
The following example shows that the converse is not true.

\begin{example}
\label{example: smooth not well formed}
Let $S=\CC[x_1,x_2,y_1,y_2,y_3]$ with
$$
\wt(x_1)=\wt(x_2)=1,\quad \wt(y_1)=\wt(y_2)=\wt(y_3)=2.
$$
Let $\PP=\Proj S=\PP(1,1,2,2,2)$.
Let $f_1$ and $f_2$ be two general polynomials of degrees $3$ and $4$, respectively.
Let $X$ be the surface in $\PP$ defined by equations $f_1=f_2=0$.
One can easily check that $X$ is quasi-smooth.
It is not well formed but weakly well formed,
because the intersection~\mbox{$X\cap \Sing \PP$} is cut out in $\PP$ by equations
\begin{equation}\label{eq:X-cap-Sing-P-11222}
x_0=x_1=f_4(0,0,y_1,y_2)=0;
\end{equation}
thus, $X\cap \Sing \PP$
is a curve on the surface $X$.
Considering the second Veronese embedding of $X$ one can check that it is smooth.
%Denote the $i$-th
%graded piece of $S$ by~$S_i$ and set
%$$
%S_{(2)}=\bigoplus_{i=0}^\infty S_{2i}.
%$$
%Let
%$$
%R=\CC[u_1,u_2,u_3,z_1,z_2,z_3]
%$$
%and $\PP'=\Proj R\cong \PP^5$.
%Consider the surjective homomorphism of the graded algebras~\mbox{$R\to S_{(2)}$} sending
%$$
%x_1^2,\ x_2^2,\ x_1x_2,\ y_1,\ y_2,\ y_3
%$$
%to
%$$
%u_1,\ u_2,\ u_3,\ z_1,\ z_2,\ z_3,
%$$
%respectively. It defines an embedding of the Veronese image $v_2(\PP)\cong\PP$
%into $\PP'$, and the image of this embedding is cut out by the equation~\mbox{$u_1u_2-u_3^2=0$}.
%Put~\mbox{$J_{(2)}=J\cap S_{(2)}$}, so that $J_{(2)}$ is the ideal of the image $X_{(2)}$ of $X$ in~$\PP'$.
%Then $J_{(2)}$ is generated by $xf_1$, $yf_1$, and $f_2$. Thus, $X_{(2)}$ is given by the equations
%\begin{multline*}
%u_1^2+u_2u_3+u_1z_1+u_3z_2=u_1u_3+u_2^2+u_3z_1+u_2z_2=\\
%\\=u_1^2+u_2^2+u_3^2+z_1^2+z_2^2+z_3^2=u_1u_2-u_3^2=0
%\end{multline*}
%in $\PP'$. One can check that these equations define a smooth variety, so that $X\subset \PP$ is smooth.
\end{example}

This example
%Example~\ref{example: smooth not well formed}
of a weakly well formed but not well formed weighted complete
intersection can be generalized to higher dimensions.

\begin{example}
\label{example: weakly not well formed}
Let $\PP=\PP(1^2,2^{N-1})$.
Let $x_1,x_2,y_1,\ldots,y_{N-1}$ be coordinates on $\PP$ with
$$
\wt(x_1)=\wt(x_2)=1,\quad \wt(y_i)=2,\ i=1,\ldots, N-1.
$$
Let $f_1$ and $f_2$ be two general polynomials of degrees $3$ and $4$, respectively.
Let $X\subset \PP$ be the weighted complete intersection defined by equations~\mbox{$f_1=f_2=0$}.
Obviously, $X$ is weakly well formed but not well formed. Indeed,
$X\cap\Sing \PP$ is defined in~$\PP$ by the equations
$$
x_0=x_1=f_2(0,0,y_1,\ldots,y_{N-1})=0.
$$
Thus, it has dimension $N-3=\dim (X)-1$ and is not a linear stratum
of $\PP$.
\end{example}

The difference between Examples~\ref{example: smooth not well formed} and~\ref{example: weakly not well formed}
for $N>4$ is that the first one is smooth, but of dimension $2$, while the second one can have any dimension, but
not smooth, and even not {quasi-smooth}.

\begin{remark}
\label{remark:generality}
The convention in~\cite{IF00} is to consider general weighted complete intersections,
that is ones given by general polynomials of given degrees. However
Examples~\ref{example: smooth not well formed} and~\ref{example: weakly not well formed} show
that well formedness and weak well formedness may differ even for complete intersections subject to the generality assumption.
\end{remark}

The notion of well formedness
%as it is introduced in Definition~\ref{definition:well formed wci}
seems to be more reasonable then weak formedness. For instance, one of the most fundamental theorem for weighted complete intersections
is the following adjunction formula.

\begin{theorem}[{see~\cite[Theorem 3.3.4]{Do82}, \cite[6.14]{IF00}}]
\label{theorem:adjunction}
Let $X\subset \PP(a_0,\ldots,a_N)$ be a quasi-smooth
well formed weighted complete intersection of hypersurfaces of degrees $d_1,\ldots,d_k$.
Then
$$
\omega_X\cong\O_X\left(\sum d_i-\sum a_j\right).
$$
\end{theorem}

In~\cite{IF00} this theorem was proved for weakly well formed varieties. However
the proof contains a gap, since in fact it uses the assumption that the variety is well formed.
For example if Theorem~\ref{theorem:adjunction} would hold for the variety $X$ from Example~\ref{example: smooth not well formed}, then one would have $-K_X\sim\mathcal O_X(1)$, so that
one gets $K_X^2=\frac{3\cdot 4}{1^2\cdot 2^3}=\frac{3}{2}$, which is impossible because $X$ is smooth.

Also, by~\cite[Proposition 2.11]{PrzyalkowskiShramov-Weighted} a smooth well formed weighted complete intersection in
a weighted projective space~$\PP$ is disjoint from the singular locus of~$\PP$; this result also fails for weakly well formed
weighted complete intersections, as shown by Example~\ref{example: smooth not well formed}.

To state the main result of this note we need the following technical definition.

\begin{definition}[{cf. \cite[Definition~6.5]{IF00}}]
\label{definition:int-with-cone}
The weighted complete intersection~\mbox{$X\subset\PP(a_0,\ldots,a_N)$}
of multidegree $(d_1,\ldots,d_k)$
is said to be \emph{an intersection
with a linear cone} if one has $d_j=a_i$ for some~$i$ and~$j$.
\end{definition}

Note that a general weighted complete intersection (cf. Remark~\ref{remark:generality}) is isomorphic to a weighted complete intersection
which is not an intersection with a linear cone, since we can isomorphically project it to
a weighted projective subspace.

It turns out that the obstructions for coincidence of notions of well formedness
and weak well formedness \emph{for weighted complete intersections} are exactly the ones that appear in Examples~\ref{example: smooth not well formed}
and~\ref{example: weakly not well formed}. More precisely, we have the following.

\begin{theorem}
\label{theorem: weakly wf is wf}
Let $X$ be a quasi-smooth weighted complete intersection of dimension at least $3$,
which is not an intersection with a linear cone. Then $X$ is well formed if and only if $X$ is weakly well formed.
\end{theorem}

The proof of this theorem is the following. By~\cite[Theorem 6.17]{IF00}, under the conditions of Theorem~\ref{theorem: weakly wf is wf}
the complete intersection is weakly well formed. Modifying the proof of loc. cit., we prove (Theorem~\ref{theorem:cone-vs-wf}) that under these conditions the complete intersection is well formed, which implies the assertion of the theorem.

Let us present some auxiliary results.

\begin{lemma}[{\cite[p. 83]{ACGH85}}]
\label{lemma: matrix of linear forms}
Let $s$, $m$, and $r$ be positive integers, and let $u$ be a non-negative integer.
Let $g_{i,j}$, $1\le i\le r$, $1\le j\le m$, be polynomials in $s$ variables, i.e. regular functions on
the affine space~$\Aff^s$.
Let $Z$ be the variety of all points $P$ in $\Aff^s$ where the matrix
$$
\left(
          \begin{array}{ccc}
            g_{1,1}(P) & \ldots & g_{1,m}(P) \\
            \vdots &  & \vdots \\
            g_{r,1}(P) & \ldots & g_{r,m}(P) \\
          \end{array}
        \right)
$$
has rank at most $u$.
Then either $Z$ is empty, or $\mathrm{codim}\,Z\le (r-u)(m-u)$.
\end{lemma}

\begin{lemma}
\label{lemma:qs-WCI-irreducible}
Let $X$ be a positive-dimensional quasi-smooth weighted complete intersection.
Then $X$ is irreducible.
\end{lemma}

\begin{proof}
Suppose that $X$ is
reducible. Let $X_1$ be one of its irreducible components.
From Hartshorne's connectedness theorem (see \cite[Theorem~18.12]{Eisenbud}) one can derive that $X$ is connected.
Let $X_2$ be another component of $X$
that intersects $X_1$ at some point $P$.
Let $C_{X_1}, C_{X_2}\subset C_X$ be affine cones over these components.
Then the intersection of $C_{X_1}$ and $C_{X_2}$ contains
the affine cone over $P$ and, thus, $C_X$ is
singular along this cone. However, this contradicts the quasi-smoothness
assumption.
\end{proof}

Now let us prove the key result of this note.

\begin{theorem}
\label{theorem:cone-vs-wf}
Let $\PP$ be a well formed weighted projective space, and let $X\subset\PP$ be a quasi-smooth weighted complete
of dimension at least~$3$
which is not an intersection with a linear cone.
Then $X$ is well formed.
\end{theorem}

\begin{proof}
Let $\PP=\PP(a_0,\ldots,a_N)$, and let $x_i$ be coordinates on $\PP$ with $\wt (x_i)=a_i$.
Let~$X$ be defined in $\PP$ by weighted homogeneous equations $f_1=\ldots=f_k=0$.
Suppose that $X$ is not well formed.
Let $\Lambda$ be a stratum of singularities of~$\PP$ such that~\mbox{$\mathrm{codim}_X (X\cap \Lambda)\le 1$}.
Without loss of generality, we may assume that $\Lambda$
is defined in $\PP$ by equations~\mbox{$x_l=\ldots=x_N=0$}.
Let $\delta$ be the greatest common divisor of
$a_0,\ldots,a_{l-1}$.
Replacing $\Lambda$ by a larger singular stratum
if necessary, we may assume that
none of the weights $a_l,\ldots, a_N$ is divisible by~$\delta$.

Suppose that $\mathrm{codim}_X (X\cap \Lambda)=0$.
Then, since $X$ is irreducible by Lemma~\ref{lemma:qs-WCI-irreducible}, we have $X\subset\Lambda$.
This means that the weighted homogeneous coordinates $x_l,\ldots,x_N$ vanish
on $X$. However, the latter is impossible because $X$ is not an intersection with a linear cone.
%(see Remark~\ref{remark:coordinate}).

Now suppose that $\mathrm{codim}_X (X\cap \Lambda)=1$.
By~\cite[Lemma 3(i)]{Di86}
one has
$$
1=\mathrm{codim}_X\left(X\cap \Lambda\right)=k(\delta)-N(\delta)+N-k+1,
$$
where $k(\delta)$ is the number of $d_j$'s divisible by $\delta$, and $N(\delta)=l$ is the number of $a_i$'s divisible by $\delta$.
This gives
$$
N-N(\delta)=k-k(\delta).
$$
Denote this non-negative integer by $r$.
Without loss of generality we may assume that $\delta$ does not divide $d_1,\ldots,d_r$.
This means that $f_1,\ldots, f_r$ vanish along $\Lambda$.
Consider the linear subspace $C_{\Lambda}\subset\Aff^{N+1}$ defined by equations~\mbox{$x_{l}=\ldots=x_N=0$};
thus,~$C_{\Lambda}$ is the affine cone of $\Lambda$.

The Jacobian matrix of the polynomials $f_1,\ldots, f_k$ at a point $P\in C_\Lambda$ has the form
$$
J(P)=
\left(
    \begin{array}{cccccc}
     0  & \ldots & 0 & g_{1,l}(P) & \ldots & g_{1,N}(P) \\
      \vdots &  & \vdots & \vdots &  & \vdots \\
     0 & \ldots & 0 & g_{r,l}(P) & \ldots & g_{r,N}(P) \\
     \vphantom{\frac{A^A}{A^A}}\frac{\partial f_{r+1}}{\partial x_0}(P) & \ldots & \ldots & \ldots & \ldots & \frac{\partial f_{r+1}}{\partial x_N}(P) \\
      \vdots &  &  &  &  & \vdots \\
     \frac{\partial f_{k}}{\partial x_0}(P) & \ldots & \ldots & \ldots & \ldots & \frac{\partial f_{k}}{\partial x_N}(P) \\
    \end{array}
  \right)
$$
for quasi-homogeneous polynomials
$$
g_{ji}=\frac{\partial f_j}{\partial x_i}\vert_{x_l=\ldots=x_N=0}
$$
in the variables $x_0,\ldots,x_{l-1}$.
Denote
$$
G(P)=
\left(
    \begin{array}{ccc}
      g_{1,l}(P) & \ldots & g_{1,N}(P) \\
      \vdots &  & \vdots \\
      g_{r,l}(P) & \ldots & g_{r,N}(P) \\
    \end{array}
  \right).
$$
Let $Z$ be the set of all points $P\in C_\Lambda\subset \Aff^{N+1}$ such that $\rk G(P)<r$. Obviously, one has~\mbox{$\rk J(P)<k$} for
$P\in Z$. The variety $Z$ is non-empty. Indeed, it contains the origin~\mbox{$0\in\Aff^{N+1}$}, because otherwise some of the polynomials
$g_{ji}$ are constant, and $X$ is an intersection with a linear cone.
Thus by Lemma~\ref{lemma: matrix of linear forms} one obtains
$$
\mathrm{codim}_{C_\Lambda}Z\le (r-(r-1))(N-l+1-(r-1))=N-N(p)-r+2=2.
$$

The affine cone $C_X$ is singular at the points of $Z\cap C_X$, that is, at the points of $Z$ subject to conditions~\mbox{$f_{r+1}=\ldots=f_k=0$}.
Denote by $Z_\PP$ the image of $Z$ in the weighted projective space $\PP$, so that $Z=C_{Z_\PP}$.
Let $S\subset X$ be the intersection of $Z_\PP$, $\Lambda$, and the common zero
set of the polynomials $f_{r+1},\ldots,f_k$ in~$\PP$.
Then
$$
\codim_\PP (S)\le
\codim_\PP Z+\codim_\PP \Lambda+(k-r)\le
2+(N-l+1)+(k-r)=k+3.
$$
Hence
$$
\dim S\ge N-(k+3)=\dim X-3\ge 0.
$$
In particular, this means that $S$ is not empty, and thus its affine cone $C_S$ is at least one-dimensional.
However, we know that the affine cone $C_X$ is singular along~$C_{S}$.
Therefore, $X$ is not quasi-smooth.
The obtained contradiction completes the proof of the theorem.
\end{proof}

\end{document}